\documentclass[12pt, a4paper]{article}
\usepackage{graphicx}
\usepackage{amssymb}
\usepackage{stmaryrd}
\usepackage{enumitem}

\usepackage{geometry}
\usepackage{tikz}
\usepackage{float}
\usepackage{url}
\usepackage{mathtools}
\usepackage{amsmath,amsfonts,amsthm}
\usepackage{mathrsfs} 
\newtheorem{definition}{Definition}
\usepackage{ytableau}
\usepackage[colorlinks,
linkcolor=blue,
anchorcolor=blue,
citecolor=blue
]{hyperref}

\begin{document}
\begin{center}
{\large \bf The Engineering and Programming Methods Used in Manufacture of Astrolabes and Errors Resulting
}
\end{center}
\begin{center}
  Duaa Abdullah\\[6pt]
  Physics and Technology School of Applied Mathematics and Informatics \\
Moscow Institute of Physics and Technology, 141701, Moscow region, Russia\\[6pt]
Email: {\tt duaa1992abdullah@gmail.com}
\end{center}
\noindent
\begin{abstract}
In this study, we first reviewed the traditional astrolabe design methods and identified potential sources of manufacturing error. We then proposed an analytical approach using computer assistance to develop designs for the astrolabe components. This approach marks a pioneering step toward designing and producing a physical astrolabe model aided by computer technology. Our goal was to revive this significant heritage instrument while leveraging modern techniques and software to produce astrolabe models free from traditional manufacturing inaccuracies. Consequently, any residual errors in the produced model are limited only by the user’s precision.
\end{abstract}


\section{Introduction}\label{sec:Intro}
The saying goes that the best way to truly understand an instrument is by making one. J.~E.~Morrison~\cite{Morrison2006} emphasize that there is very limited information on how to craft an astrolabe. The expertise and tools needed to create a functional metal astrolabe involve a complex science that exceeds the scope of this work. Moreover, techniques for working with and engraving metal instruments are continually evolving, making it impractical to cover a significant portion of these methods in meaningful detail here.
The \emph{astrolabe’s} flexibility and user-friendliness stem from the method it employs to depict the sky on its face. Visualizing the sky is inherently challenging. Throughout history, people have \emph{gazed upward}, striving to comprehend its complexities. Since our viewpoint is restricted to a \emph{single location} and limited by changing conditions, we never see more than a portion of the sky at any moment. The visible features are continuously shifting: the Sun and stars rise and set each day, the array of visible stars changes with the seasons, and the Sun’s position in the sky varies from day to day. Throughout the centuries, a variety of instruments have been employed for both orientation and measurement, but the astrolabe stands out as the one most capable of fulfilling these functions with precision. The astrolabe represents a flat projection of an armillary sphere; it models the apparent rotation of the celestial sphere relative to the Earth and a specific latitude, allowing the user to determine the positions of stars at any given time. A detailed description of its construction and functionality is deliberately omitted here, as it is covered comprehensively in a prior study on the subject~\cite{Aterini2019}. \par
\begin{figure}[H]
    \centering
    \includegraphics[width=0.5\linewidth]{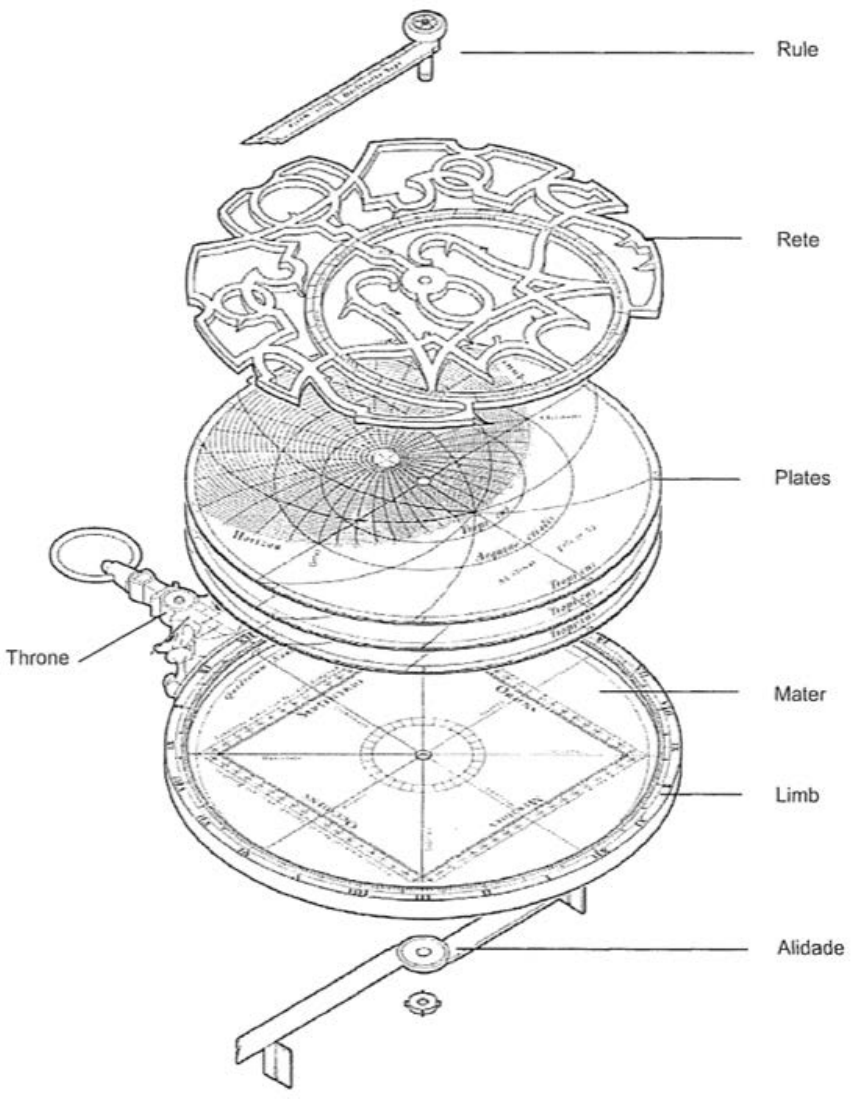}
    \caption{Parts of Astorlabe.}
    \label{figast01}
\end{figure}

 
The largest circle on Earth that lies perpendicular to the planet’s axis is called the equator; all smaller circles that run parallel to the equator are known as \emph{parallels}, Great circles that pass through both poles are called \emph{meridians}, with the one passing through Greenwich designated as the prime meridian $0^{\circ}$.
The flat astrolabe is one of the greatest achievements of Arab and Muslim scientists; it is a computational observational instrument. Despite its presence in museums, it has not lost its scientific importance, as it remains usable even today. It serves as an observational device \emph{theodolite}, performs calculations \emph{calculator}, and can even be said to solve certain problems like programmed calculators. Another aspect is the ability to solve problems directly and obtain results from the astrolabe itself without needing to write anything on paper. Figure~\ref{figast01} shows its parts. The astrolabe can be used as an observational instrument to determine altitudes and deviations, and as a computational device capable of solving hundreds of astronomical and surveying problems related to observation. For this reason, it is called the noble instrument.

The assumption that readers already understand astrolabe principles and the different scales is common in many scholarly papers, but it may not hold true even for advanced audiences. Including every variation of every scale ever used on astrolabes in a reasonably sized book is impractical. Nevertheless, the content presented offers a solid foundation for interpreting historical instruments and can serve as a useful starting point for deciphering the purpose and design of scales not explicitly discussed. 
The design of the astrolabe relies on the theory of sphere flattening (projection), which was developed by Arab and Muslim scientists. Consequently, complex calculations in space (on the sphere's surface) were transferred to the plane, representing a creative work and a great achievement that manifested in simplifying the solution of complex astronomical problems.

The following symbols are used throughout this paper:
\begin{center}
\begin{tabular}{|l|l|}
\hline
\textbf{Symbol} & \textbf{Definition} \\ \hline
$h$, $a$     & altitude of a celestial body above the horizon \\ \hline
$A$          & azimuth of a celestial body \\ \hline
$H$          & hour angle of the Sun or star \\ \hline
$\delta$     & Declination \\ \hline
$\alpha$     & Right ascension \\ \hline
$\varphi$    & terrestrial latitude \\ \hline
$\beta$      & celestial latitude \\ \hline
$\lambda$    & celestial or terrestrial longitude \\ \hline
$\varepsilon$& obliquity of the ecliptic \\\hline
\end{tabular}
\end{center}

This paper prepared as follows. In section~\ref{sec:Projection}, we presented the projection system are used in  design the astrolabe. Section~\ref{sec:Traditional} discused the popular method are used in design astorlabe mainly use geometric construction and stereographic projection. Section~\ref{sec:Stereographic} discussed stereographic design using computer assistance.

\section{The Projection System Used in the Astrolabe}~\label{sec:Projection}
The astrolabe is considered one of the most important astronomical instruments that contributed to the development of astronomy, as it is regarded as the foundation upon which astronomy was built and served as a kind of computer that made life easier at that time, which was the golden age of the astrolabe. Due to its importance, it was necessary for us to clarify the methods of manufacturing the astrolabe; however, because of the many types of astrolabes, we have chosen the flat astrolabe in this paper.\par
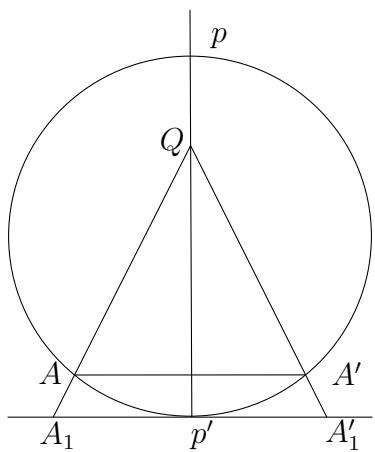
\begin{figure}[H]
    \centering
  \begin{tikzpicture}[scale=.6]
\draw   (3,1)-- (11,1);
\draw   (7,5) circle (3.9781905434506273cm);
\draw   (7,10)-- (7.04,1);
\draw   (7.013333069964051,7.000059258088729)-- (4,1);
\draw   (7.013333069964051,7.000059258088729)-- (10,1);
\draw   (4.467910885922933,1.9316902505821867)-- (9.535023337064967,1.9341140463911595);
\draw (7.22,9.84) node[anchor=north west] {$p$};
\draw (6.1,7.62) node[anchor=north west] {$Q$};
\draw (3.44,2.46) node[anchor=north west] {$A$};
\draw (9.86,2.46) node[anchor=north west] {$A^{\prime}$};
\draw (3.46,1.12) node[anchor=north west] {$A_1$};
\draw (9.72,1.2) node[anchor=north west] {$A_1^{\prime}$};
\draw (6.78,1.18) node[anchor=north west] {$p^{\prime}$};
\end{tikzpicture}
    \caption{Projection point of Astorlabe.}
    \label{figast02}
\end{figure}

The astrolabe employs a projection method known as stereographic projection. This approach transforms the three-dimensional celestial sphere onto a flat, two-dimensional surface, usually by projecting from the South Celestial Pole onto the sphere’s equatorial plane. The tympan (plate) of the astrolabe features this projection, engraved with altitude (almucantar) circles, azimuth lines (meridians), the horizon, zenith, and additional reference circles essential for astronomical calculations and navigation.
Perspective projection is considered one of the primary methods in map drawing. It was used to directly project the surface of the celestial sphere onto a plane from a viewpoint (projection point) $Q$ located on the celestial axis $pp^{\prime}$ or its extension, as shown in Figure~\ref{figast02}.
will start with standard points among Definition~\ref{deffstan1} and lines and continue on with details covering each of the additional four map projection parameters.
\begin{definition}[Standard point~\cite{Ingram}]~\label{deffstan1}
A standard point and line is a point or line of intersection between the developable surface and the spheroid or ellipsoid. In the case of a secant intersection, there will be two standard lines that would define where the developable surface intersects with the spheroid.
\end{definition}
Depending on the distance from point $Q$ to the projection plane, four types of projection methods are defined:

\begin{itemize}
    \item \textbf{Orthographic Projection\footnote{Furthermore, see \href{https://uw.pressbooks.pub/enggraphics/chapter/orthographic-projection/}{Orthographic Projection}.}:}
     is a technique for depicting three-dimensional objects on two-dimensional surfaces by projecting the object's features onto perpendicular planes using parallel lines that intersect these planes at right angles. This creates several views of the object—commonly three: front, top, and side—each displaying the object's dimensions accurately without perspective distortion. This method is commonly used in engineering and architecture because it maintains the true size and shape of the object's surfaces, enabling precise measurement and construction. The projection point $Q$ in Figure~\ref{figast02} is very far from the projection plane, making the projection lines parallel. This is the basis for imaging distant objects, such as representing the surfaces of planets and the Moon. 
    
    \item \textbf{External projection}: External projection in cartography is a map projection method where the viewpoint or projection center is positioned outside the Earth. This approach creates a perspective that mimics observing the Earth from a point in space, some distance away from the planet~\cite{Ingram}. The projection point $Q$ is located outside the sphere at a specific distance. This is the basis for monoscopic space images from satellites.
    
    \item \textbf{Stereographic projection}: In the stereographic projection~\cite{Ingram}, the light sources are positioned on the side of the Earth opposite to where the developable surface intersects secantly or tangentially. This setup results in a milder distortion between the compressed and stretched areas of the Earth, although no region is entirely free of distortion. The projection point $Q$ coincides with the celestial pole $p$. This projection is used in creating maps and in designing the astrolabe, the subject of our research.
    
    \item \textbf{Gnomonic projection}: This is the case where the projection point $Q$ is at the center of the sphere. This projection was used for some astronomical maps~\cite{Chau1987Dean}.
\end{itemize}

After presenting the general principle of perspective projection methods, we pose the following question: Why was stereographic projection chosen as the principle in designing the planispheric astrolabe?

In fact, attempting to answer this question reveals the vast scientific imagination of Arab and Muslim scholars. 
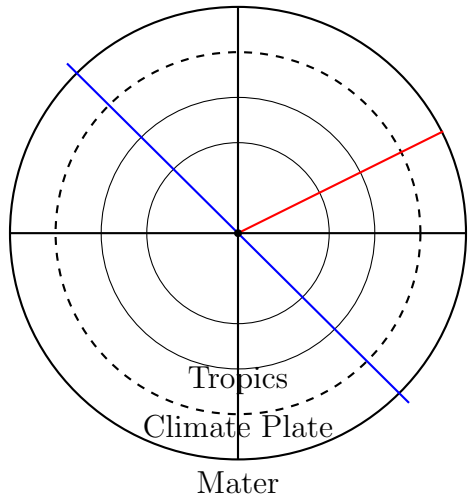
\begin{figure}[H]
\centering
 \begin{tikzpicture}[scale=1.5]
\draw [thick] (0,0) circle (2);
\draw [thick, dashed] (0,0) circle (1.6);
\draw [thin] (0,0) circle (1.2);
\draw [thin] (0,0) circle (0.8);
\draw [thick] (0,2) -- (0,-2);
\draw [thick] (-2,0) -- (2,0);
\draw [red, thick] (0,0) -- (1.8,0.9);
\draw [blue, thick] (-1.5,1.5) -- (1.5,-1.5);
\node at (0,-2.2) {Mater};
\node at (0,-1.7) {Climate Plate};
\node at (0,-1.3) {Tropics};
\filldraw (0,0) circle (0.03);
\end{tikzpicture}
\caption{Circles and straight lines of Astorlabe.}
    \label{figast03}
\end{figure}
They realized that in this projection, the following property is achieved: \emph{the projection of any circular arc is also a circular arc}. 
Consequently, lines on the surface of the sphere that have astronomical significance—such as the equator (the ecliptic of Aries), the tropics of Cancer and Capricorn, lines of longitude (hour circles), circles parallel to the horizon (almucantars), azimuth lines (altitude circles), the zodiac circle, and the orbits of other planets—have projections that are all circles. 
This property offers immense benefit, manifested in the ease and accuracy of drawing these arcs when manufacturing the astrolabe plate.

To emphasize the importance of this property, let us consider the case where the projection point deviates from the pole while remaining on the extension of the axis. 
In this case, the stereographic property disappears. 
The projections of inclined astronomical circles will no longer be circular arcs but will become conic sections~\cite{Jasem3}, namely ellipses, hyperbolas, and parabolas, in addition to circles and straight lines. Figure~\ref{figast03} illustrates one of these cases when the projection point is at the center of the sphere. Were Arab and Muslim scholars able to envision these projections about nine centuries ago to choose the correct or appropriate location for the projection point? Al-Biruni\footnote{Al-Biruni had some ideas very strikingly similar to those of Einstein and other modern scientists regarding the Universe as a whole~\cite{Ahmad2010}.} says in this regard:
``And if the position of the operation is changed, meaning the pole of the plane, without deviating from the alignment of the axis, the divisions then change and become straight lines and circles (see Figure~\ref{figast03}) and ellipses, hyperbolas, and parabolas.''

Thus, the stereographic projection~\cite{Howarth1996,Lisle2004} is considered the only projection among all known projections, ancient and modern, that possesses this property. 
There is no other projection that ensures the transfer of any circular arc from the surface of the sphere to the projection plane in the form of a circular arc, except for this projection.

\section{Traditional Methods in Designing and Errors}~\label{sec:Traditional}
Popular methods for designing astrolabes mainly use geometric construction and stereographic projection. Traditional design involves precise drafting guided by mathematical principles, where the celestial sphere is projected stereographically onto a flat plate called the tympanum, tailored to a specific latitude. The process includes drawing perpendicular lines and circles with a compass and straightedge to achieve accurate calibration for measuring celestial positions. An astrolabe's main components are the mater (the base), tympans (plates customized for different latitudes), and the rete (a rotating star map), all crafted using these geometric and mathematical techniques. The design aims to mimic the apparent motion of stars and the sun for the user's latitude, allowing diverse astronomical calculations.
In general, solving any astronomical problem using the astrolabe requires observation and calculation. 
Observation involves obtaining measurements using the alidade and the back of the astrolabe. 
Calculation, on the other hand, involves deriving other measurements or results based on the observation data. 
Therefore, the solution to any problem must include a composite error from two sources: the observation process and the calculation process. 
In line with this, astrolabe errors can be classified into two categories:
\begin{itemize}
    \item \textbf{The first:} Errors in the calculation elements, represented by the face of the astrolabe.  
     \item \textbf{The second:} Errors in the observation elements, generally represented by the back of the astrolabe.
\end{itemize}

\subsection{First: Potential Errors in the Calculation Elements}
\begin{itemize}
    \item \textbf{Plate errors:}
    \begin{itemize}
        \item Radii of the Tropics of Cancer, Capricorn, and the Equator.
        \item Corresponding almucantars and their graduations.
        \item Altitude circles.
        \item Hour circles.
    \end{itemize}
    \item Eccentricity error of the mater and its graduations.
    \item \textbf{Rete (spider) errors:}
    \begin{itemize}
        \item Design of the ecliptic circle and its graduation.
        \item Star pointer fragments.
    \end{itemize}
\end{itemize}

\subsection{Second: Potential Errors in the Observation Elements}
\begin{itemize}
    \item Unequal squares.
    \item Unequal altitude divisions (graduations of vertical angles).
    \item Non-horizontality of the reference line for measuring altitude angles (the east-west line).
    \item Error in the shadow square.
    \item \textbf{Alidade errors:}
    \begin{itemize}
        \item Misalignment of the sighting axis (axis of the sight vanes) with the line of the pointers.
        \item Misalignment of the alidade's rotation axis with the line of the pointers.
    \end{itemize}
\end{itemize}
We will address these errors through geometric analysis.

\subsection{Plate Errors}
\begin{figure}[H]
\centering
    \includegraphics[width=0.9\textwidth]{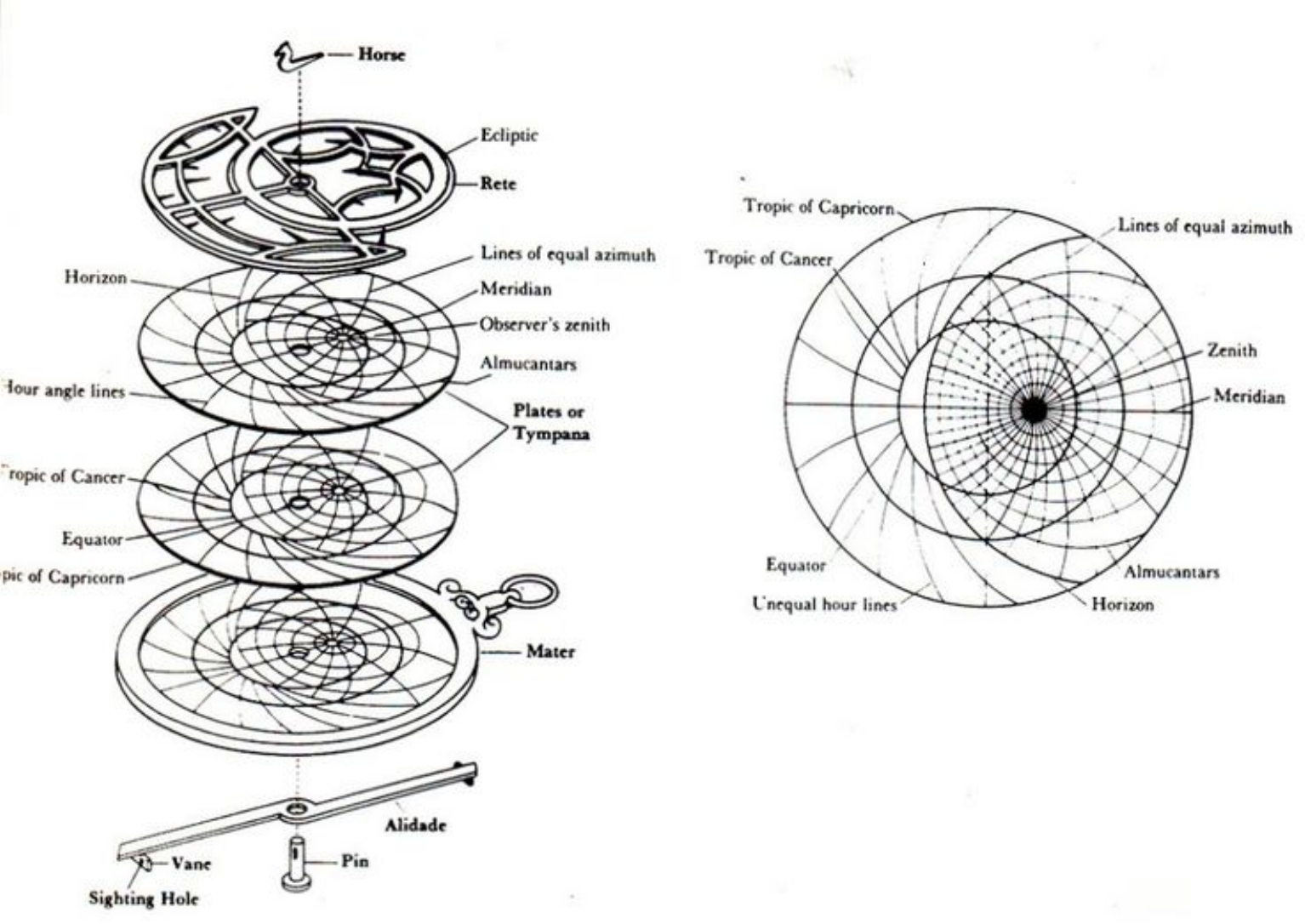}
 \caption{Almucantars, azimuths, and hours of astorlabe.}
\label{figp01}
\end{figure}

The plate is considered the most important part of the astrolabe; it is a disk on which curves called almucantars, azimuths, and hours are drawn, as shown in Figure~\ref{figp01}. Plate errors in astrolabes primarily stem from inaccuracies in scale graduations, calibration challenges, and observational constraints. The plate’s divisions often have slight positional errors—usually a few arc minutes—caused by irregular engraving and the physical width of the markings. These inaccuracies can increase depending on the user’s ability to precisely align the alidade or sighting device when measuring the altitude of celestial objects. Additionally, some systematic errors arise from the plate’s projection method, which is tailored for a specific latitude and can introduce distortions if used outside that region. While finely made astrolabes can achieve accuracy within a few arc minutes, typical practical errors frequently range from 10 to 15 arc minutes due to these combined factors.
The design principle of the plate relies on the stereographic projection of a spherical coordinate system onto a plane. The first coordinate system is the set of celestial latitude and longitude lines that take the North Pole as their origin, and it appears on the plate through three concentric circles centered at the plate's center, representing the tropics of Capricorn, Cancer, and Aries (the equator).

The second coordinate system is the set of altitude and azimuth circles that take the zenith as their origin, and it appears on the plate through curves called almucantars and azimuths. The type of projection (stereography) used determines the nature and shape of the curves on the plate.

\begin{itemize}
    \item \textbf{Errors in the Radii of Aries, Cancer, and Capricorn\footnote{At present, these orbits cannot be called Aries and Cancer due to the deviation of the zodiac circle. The orbits of Capricorn, Aries and Cancer have become Sagittarius, Pisces and Gemini, respectively.}:}
Since the radii of the aforementioned circles are determined graphically, there are undoubtedly errors committed in drawing them. To detect these errors, consider Figure~\ref{figp02} and the projections of these tropics in it. We observe that if the sphere's diameter is taken as one unit, then according to stereographic projection, the radius of the Capricorn tropic in the projection should be $1.53$, the Aries tropic $1$, and the Cancer tropic $0.66$. Thus, the distance between Capricorn and Aries is $0.53$, and between Capricorn and Cancer is $0.87$.

\begin{figure}[H]
\centering
\begin{tikzpicture}[scale=.5]
\draw  (5.02999024390244,0) circle (3cm);
\draw [color=red] (5.02999024390244,0) circle (4cm);
\draw [color=blue] (5.02999024390244,0) circle (5cm);
\draw   (18,6)-- (18,-6);
\draw   (20,0) circle (2cm);
\draw   (20.838449008719028,1.8157651995172948)-- (20.828305213857327,-1.820414917730782);
\draw   (20,2)-- (20,-2);
\draw   (18.876231720117218,1.6544318822874797)-- (18.822442938057296,-1.6165888054384354);
\draw   (5.021218011846493,2.999987174630044)-- (18,2.8670048780487805);
\draw   (18,2.8670048780487805)-- (22,0);
\draw   (5.02999024390244,4)-- (18,4);
\draw   (18,4)-- (22,0);
\draw   (4.99533778388403,4.999879919259529)-- (18,4.814809756097562);
\draw   (18,4.814809756097562)-- (22,0);
\draw [->] (6.93504957089905,-2.3174876397994497) -- (9.081424390243905,-4.067180487804882);
\draw [->] (6.935220671901679,-3.517114871059209) -- (9.445014634146345,-5.703336585365858);
\draw [->] (6.7067173634910375,-4.710476214400846) -- (9.57486829268293,-7.365463414634152);
\draw (9.13336585365854,-3.391941463414638) node[anchor=north west] {$Cancer$};
\draw (9.445014634146345,-5.028097560975614) node[anchor=north west] {$Equator $};
\draw (9.8345756097561,-6.8) node[anchor=north west] {$Capricorn$};
\end{tikzpicture}
\caption{Orbits in the zodiac.}
\label{figp02}
\end{figure}
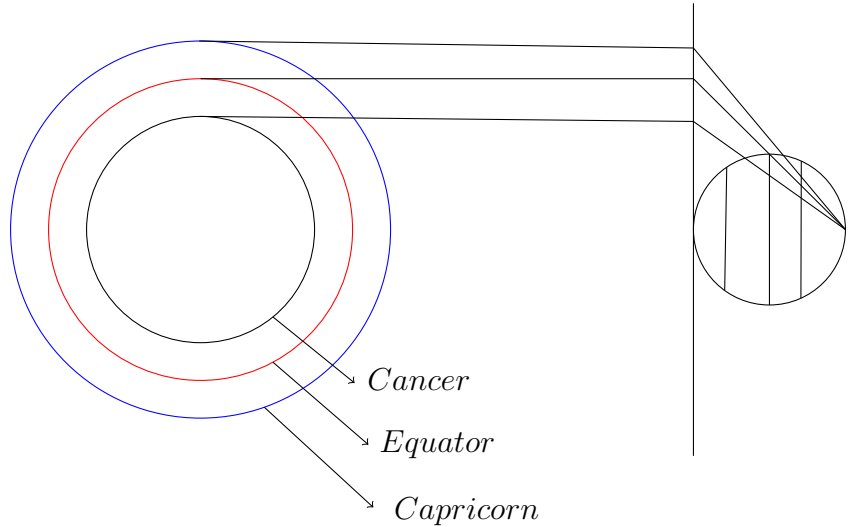
\item \textbf{Almucantars}: An almucantar is a small circle on the celestial sphere that runs parallel to the horizon, connecting all points sharing the same altitude. Simply put, any two stars or celestial objects on the same almucantar have identical altitudes above the horizon. The term originates from the Arabic word "al-muqanṭarah," meaning arch or bridge. Almucantars play a crucial role in astronomical instruments like the astrolabe, where they represent altitude circles used to measure the height of celestial bodies relative to the observer’s horizon. These circles appear as concentric rings on the astrolabe’s tympanum, enabling accurate determination of celestial altitudes and supporting precise astronomical observations and calculations (see Figure~\ref{fig003Tierra}).
\end{itemize}
Al-Marrakushi\footnote{\textbf{Ibn al-Bann\={a}' al-Marr\={a}kush\={\i} (1256--1321)} was a Moroccan polymath renowned as a mathematician, astronomer, Islamic scholar, Sufi, and occasional astrologer. Born in Marrakesh, he contributed significantly to arithmetic, algebra, and astronomy during the Islamic Golden Age, including rediscovering Th\={a}bit pairs and authoring texts like \textit{Talkh\={\i}\d{s} \={a}\d{m}al al-\d{h}is\={a}b} (on fractions and sums of powers) and \textit{Raf\d{} al-\d{h}ij\={a}b} (on square roots and continued fractions, featuring early algebraic notation). His works influenced later scholars and North African mathematics.} mentioned in his explanation of the method for drawing almucantars that the horizon must first be determined. If we assume the observation point $S$ has geographic latitude $\phi$ (Figure~\ref{figp03}), then the zenith $Z$ lies at an angular distance of $(90^\circ - \phi)$ from the pole $P$.

\begin{figure}[H]
\centering
\includegraphics[width=0.7\textwidth]{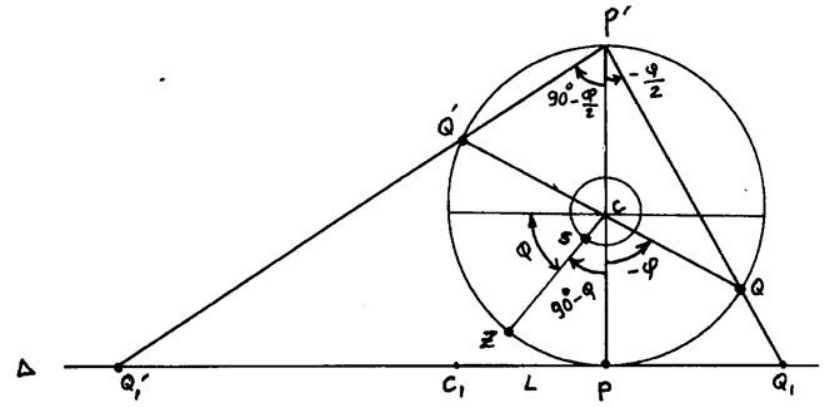}
\caption{Almucantars.}
\label{figp03}
\end{figure}

\begin{figure}[H]
\centering
\includegraphics[width=0.7\textwidth]{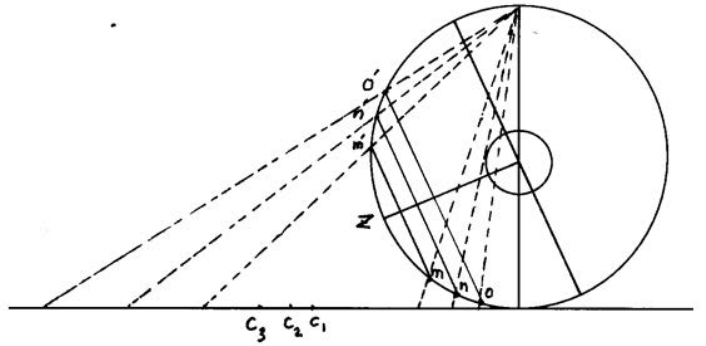} 
\caption{Draw almucantars.}
\label{figp04}
\end{figure}

To draw the projections of the almucantars, we first define the almucantar on the sphere by taking points at specific intervals (for example, every ten degrees) starting from $Z$ (or from the horizon) and symmetrically in both directions (Figure~\ref{figp04}). Then, we connect each pair of symmetric points to obtain the almucantars on the sphere. To find their projection, we follow the same method as projecting $QQ^{\prime}$, yielding their radii and centers.

Hence, we observe the accumulation of errors during the drawing of the almucantars.

\subsection{Altitude Circles}

According to al-Marrakushi's description (Figure~\ref{figp03}), the angle $\angle BOB_3$ equals the angle $\angle B_1OB_2$ by perpendicularity, which indicates the geographical latitude $\varphi$ of the location\footnote{Altitude circles (also called \emph{circles of equal altitude}) are small circles on the Earth's surface where a celestial body has the same altitude (angular height above the horizon) at a given time.}. 
\begin{figure}[H]
\centering
  \begin{tikzpicture}[scale=.5]
\draw   (7,3) circle (3cm);
\draw   (9.019013168062262,5.218915461929812)-- (5.061248864012673,0.7106236585681573);
\draw   (4.878679656440358,5.121320343559643)-- (9.121320343559642,0.8786796564403576);
\draw   (7,6)-- (7,0);
\draw   (4,3)-- (10,3);
\draw   (-4,0)-- (7.026949807077819,5.999878948874182);
\draw   (7.026949807077819,5.999878948874182)-- (12,0);
\draw   (12,0)-- (-4,0);
\draw (-4.302,0.0368) node[anchor=north west] {$T$};
\draw (11.7426,-0.0842) node[anchor=north west] {$T_1$};
\draw (6.6364,0.1578) node[anchor=north west] {$A$};
\draw (-0.309,0.0126) node[anchor=north west] {$A_1$};
\draw   (0,0)-- (0,6);
\draw (0.0782,6.3046) node[anchor=north west] {$A_2$};
\draw (6.7574,7.1032) node[anchor=north west] {$B$};
\draw (9.008,6.111) node[anchor=north west] {$B_1$};
\draw (10.1938,3.7878) node[anchor=north west] {$B_2$};
\draw (4.0712,6.4) node[anchor=north west] {$B_3$};
\draw (6.146,4.5138) node[anchor=north west] {$\varphi$};
\draw (7.6528,4.0056) node[anchor=north west] {$\varphi$};
\end{tikzpicture}
  \caption{Almucantars provided in Figure~\ref{figp03}}
\label{figp06}
\end{figure}
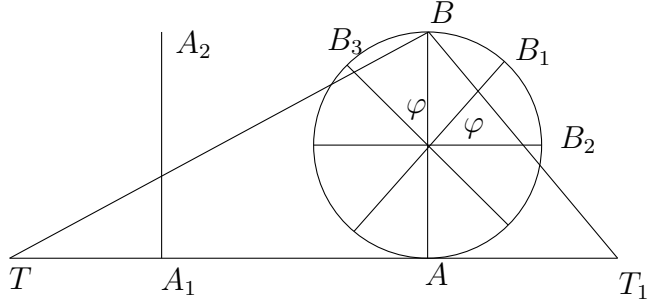
Then, $B_1B_2$ becomes the zenith, as it is the midpoint of the arc $BB_1B_2$, because $B_1B_2$ is the horizon of the location. 
As for $A_1$, it is the nadir point. 
The altitude circles are the great circles passing through the two points $k$ and $A_1$. 
To find the projection of these circles, we project the zenith point $B_1B_2$ and the nadir $A_1$ using stereographic projection from $B$, obtaining the points $T_1$ and $T$.  Consequently, the projections of the altitude circles must pass through the projections of the zenith and nadir points, because the stereographic projection process is, in reality, an optical projection operation. 

Since $TT_1$ represents the diameter of the circle $B_3BB_2$ in the projection, its center $O$ can be easily found by halving $TT_1$. The circle $kmd$ intersects the horizon circle and the equator circle $B_3BB_2$ on the sphere (Figure~\ref{figp06}).

 This intersection will appear in the projection at the points $A$ and $T_1$. The circle $AT_1$ is called the prime altitude circle. Thus, we now have circles passing through the two points $TT_1$, one of which has its center at $O$ as mentioned earlier. As for the remaining circles, $TT_1$ represents a chord in them. Therefore, the centers of these circles lie on the axis of the segment $TT_1$, that is, on the perpendicular erected on $TT_1$ at $O$.

Now, we must draw circles through $TT_1$ such that their centers lie on $O$ with a specific angular interval (ten degrees, for example) from the prime altitude circle. 

\subsection{Hour circles}

\begin{figure}[H]
    \centering
    \includegraphics[width=0.5\linewidth]{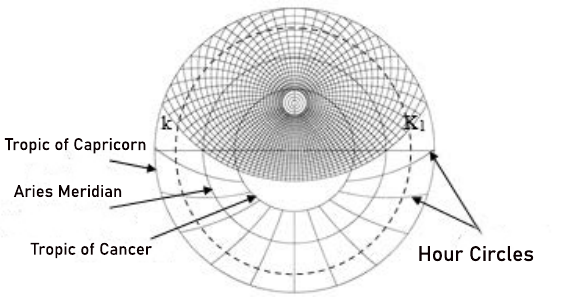}
    \caption{Hour Circles}
    \label{fig001HourCircles}
\end{figure}

The principle of dividing the hour arcs is based on dividing the orbits of Cancer, Capricorn and the equator into twelve equal parts, then connecting these parts with circular arcs. 
However, the accuracy of this construction (according to stereographic projection) requires that the division arc also be circular on the surface of the sphere. 
In reality, however, the arcs of division on the surface of the sphere are not flat but left-handed, and therefore their projection is not circular.
However, the error resulting from considering them circular can be neglected, as it is often less than the accuracy of the astrolabe. 
For the astrolabe maker, drawing the hour arc means passing a circular arc through three given points, i.e. choosing an arc with a suitable radius and centre, which is difficult without the use of a computer.
Testing the accuracy of the division of the hours depends on matching the movement of one part of the zodiac by one hour with the graduations obtained from the parts of the chamber, for each of the twelve hours. 
This means that any part of the zodiac must cut equal angles between the limits of each hour. 
Therefore, each arc $kk_1$ between the tropics of Cancer and Capricorn (Figure~\ref{fig001HourCircles}) must be divided into twelve equal parts. 
If there is an error in the division, its effect will appear in solving time-related problems. 
For example, in determining the time remaining until sunset, we observe the altitude of the sun and determine its position on the altitude circles, and two readings of time are taken, the first at this altitude and the second when the sun crosses the horizon circle. 
Here, the error in the answer to the problem comes from the error in the altitude circle (arch), the error in the horizon circle, and the error in the hour arcs. Note that there are multiple sources of error. 
If we consider them to be random errors, they may accumulate in a consistent manner, leading to an incorrect result.

\subsection{Error in the decentralisation of the chamber and its gradations}
\begin{figure}[H]
    \centering
    \includegraphics[width=0.4\linewidth]{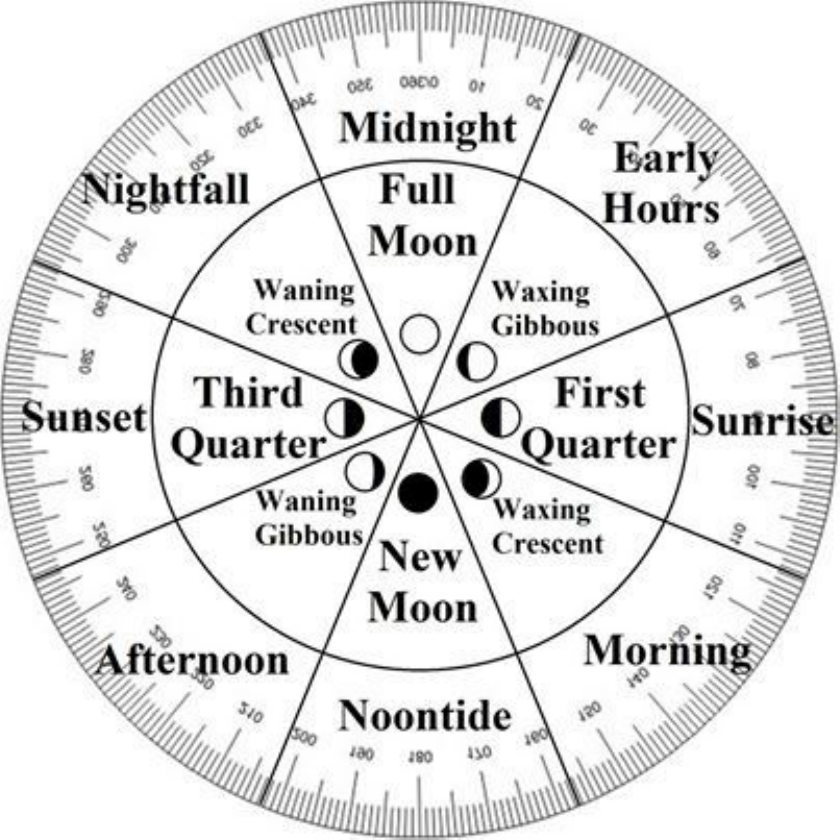}
    \caption{A fundamental phases  in Astorlabe.}
    \label{fig002fundamental}
\end{figure}

\begin{itemize}
    \item \textbf{Unequal quadrants:}
\end{itemize}
The fundamental error in the design and construction of the chamber is the inequality of the quadrants, which results either from the vertical axis (midday line) not being perpendicular to the line of elevation measurement (sunrise and sunset line), or from the decentralisation of the circle of gradations relative to the chamber. 
This error can be detected by measuring the chords of these quadrants with a fixed compass, as shown in Figure~\ref{fig002fundamental}. 
If the adjacent quadrants differ while the opposite ones are equal, this indicates that the sunrise-sunset line is not horizontal. 
If all the quarters differ, this reveals the eccentricity of the circles of gradations (the chamber and the back of the astrolabe) with respect to the polar axis. 

\begin{itemize}
    \item \textbf{Gradations:}
\end{itemize}
 
In addition to the quarter error, there is also a gradation homogeneity error. This error can be constant or variable, and can be detected based on the principle that equal chords on the circumference of a circle correspond to equal arcs and equal central angles. 

\subsection{Network (spider) errors}
we will discuss the following issues:
\begin{enumerate}
    \item[1] \textbf{The zodiac circle and its gradations:}
\end{enumerate}
The zodiac division test is based on the occurrence of a constellation and its counterpart on the ascending and descending arcs, or on the ascending and descending lines, or on the midday line. 
Considering the equal division of the zodiac circle on the celestial sphere, each part of the zodiac and its counterpart are located at opposite ends of a great circle, inclined at an angle of approximately 23.5 degrees to the equator. 
Considering that the horizon is also a great circle, at any given moment there are six zodiac signs above the horizon and six below it. 
This means that if any part of the zodiac is on the eastern horizon, its counterpart must be on the western horizon. 
On the other hand, since each part of the zodiac and its counterpart lie at opposite ends of a great circle, this circle must coincide with the lines of sunrise and sunset and the meridian. 
Therefore, when any part of the zodiac falls on the line of the east, its counterpart must fall on the edge of the west. 
If this does not happen, it means that there is an error in the divisions of the zodiac circle.

\begin{figure}[H]
    \centering
    \includegraphics[width=0.5\linewidth]{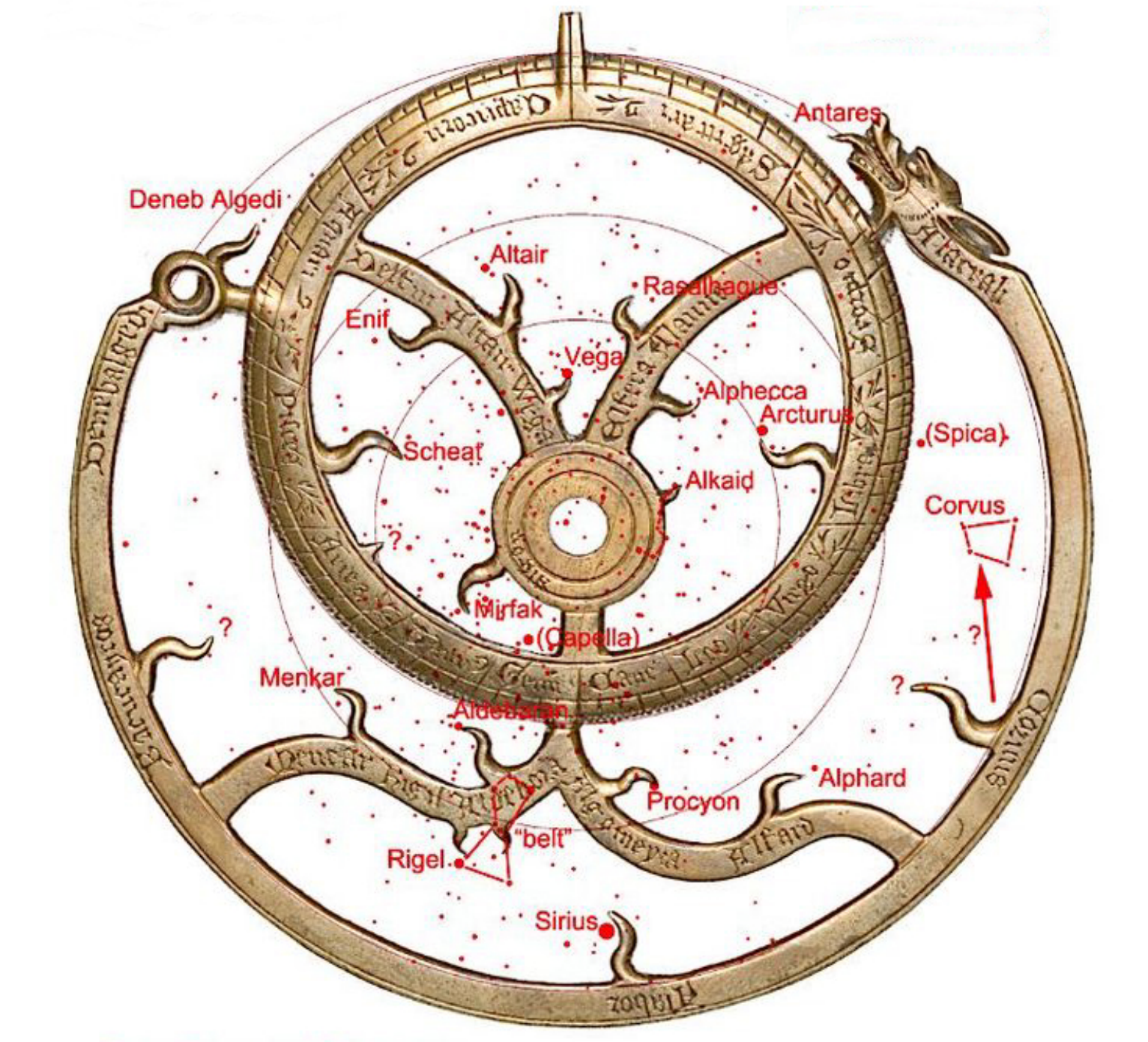}
    \caption{Tierra hueca in Astorlabe.}
    \label{fig003Tierra}
\end{figure}

\begin{enumerate}
    \item[2] \textbf{Fragments of fixed stars:}
\end{enumerate}
Errors in the fragments of stars mean errors in their coordinates (Figure~\ref{fig003Tierra}). Therefore, this error is detected in the following sequence:
\begin{itemize}
    \item  Measure the actual altitude of one of the stars.
    \item Place its fragment on the arc of this altitude on the astrolabe.
    \item Read the altitude of another star at the same moment from the plate.
    \item Quickly observe the latter in reality.
    \item Read its altitude from the gnomon.
    \item Compare this reading with the previous reading from the plate; the two readings must match. This match means that the relative position between the two stars is correct, i.e. their coordinates are correct. Repeat these steps for the rest of the fragments to verify the coordinates of the remaining stars. 
\end{itemize}

\begin{enumerate}
    \item[3] \textbf{The astrolabe and its manufacturing errors:}
    \begin{enumerate}
        \item Unequal quadrants,
        \item Unequal elevation segments (vertical angle increments),
        \item Non-horizontal line of elevation angle measurement (sunrise and sunset line). 
        \item Error in the shadow square: The accuracy of the shadow square is very important, because most surveying problems depend on observation with a staff and the use of a shadow square. In this box, the divisions must correspond to equal angles, especially the diameter of the box with the line passing through the 45-degree elevation scale. Considering that the shadow function is not linear, its divisions will not be homogeneous, meaning that errors will inevitably occur in its division using traditional methods. 
        \item Errors in the hinge (Figure~\ref{fig004rate}):
        \begin{enumerate}
         \item The axis of rotation of the hinge does not coincide with the line of the two fragments.
         \begin{figure}[H]
             \centering
             \includegraphics[width=0.7\linewidth]{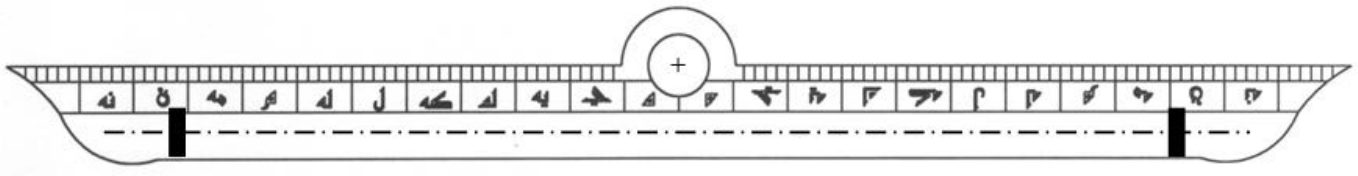}
             \caption{Rate in Astorlabe.}
             \label{fig004rate}
         \end{figure}
         When designing and manufacturing the spline, it is possible to make a mistake in determining its centre of rotation. This mistake is a shift in a certain direction from the correct centre, and can be analysed into two perpendicular components, one of which applies to the line of the splines and the other to the midday line. Therefore, the existence of this eccentricity can be tested in two stages. We apply fragment A to the sunset line (Figure 16), and fragment B must apply to the sunrise line, otherwise this reveals an eccentricity towards the midday line of 
         \begin{equation}~\label{eqq1rate}
             d_1=\frac{\ell\delta}{4},
         \end{equation}
         where $\ell$ is the length of rate and $\delta$ is the difference between the head and line orientation. In Figure~\ref{fig005Error}, we observe the error while rotate the rate by considering to $\epsilon$ where $\epsilon$ is the difference between other lines according to~\eqref{eqq1rate} such that
         \begin{equation}~\label{eqq2rate}
             d_2=\frac{\epsilon}{4}.
         \end{equation}
         \begin{figure}[H]
             \centering
             \includegraphics[width=0.5\linewidth]{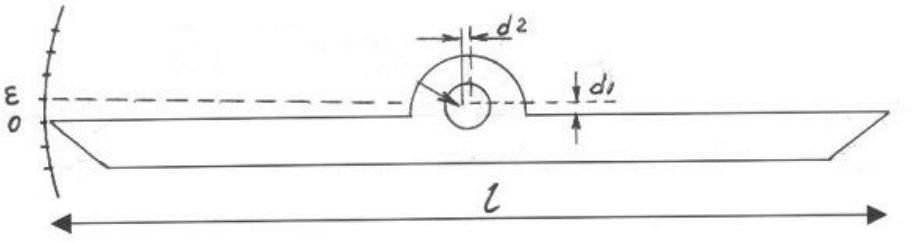}
             \caption{Error in Rate.}
             \label{fig005Error}
         \end{figure}
\item The observation axis does not align (or parallel) with the line of the fragments:
Testing for the presence of an error $\alpha$  in the observation axis  - which is essential in the manufacture of the crossbar –  is done by observing the height of a point at both ends of the prism and at both ends of the height division circle. The two readings must match, otherwise this reveals the presence of error $\alpha$  between the axis of the two targets and the line of the prisms.
\end{enumerate}
\end{enumerate} 
\end{enumerate}

\section{Analytical Study of Systematic Errors in Drawing Arcs in the Planimeter}

\begin{figure}[H]
    \centering
    \includegraphics[width=0.4\linewidth]{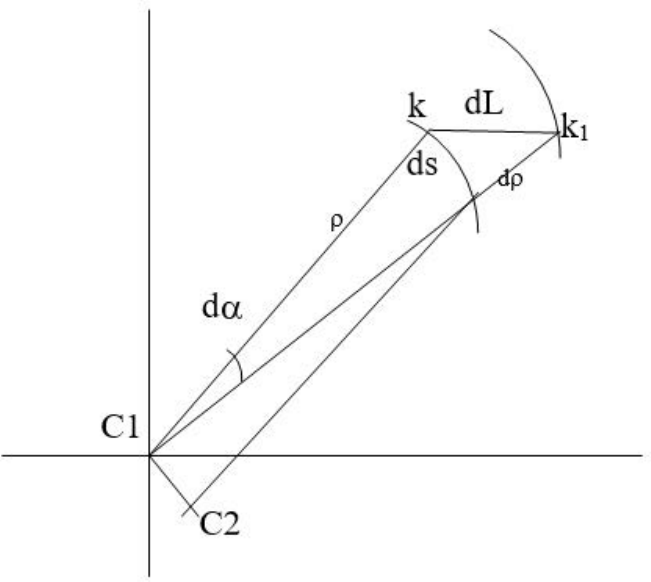}
    \caption{Drawing arcs in the planimeter.}
    \label{fig006Drawing}
\end{figure}

The deviation of any point $k$ from a circular arc from its correct position in the plane to $k_1$ is due to an error in determining the center of the arc from $C_1$ to $C_2$, and for a radius increasing from $\rho$ to $\rho + dp$, as shown in Figure~\ref{fig006Drawing}. The displacement can be calculated~\eqref{eqq3rate} as follows:
\begin{equation}~\label{eqq3rate}
dL = \sqrt{ds^2 + dp^2}.
\end{equation}

\begin{figure}[H]
    \centering
    \includegraphics[width=0.4\linewidth]{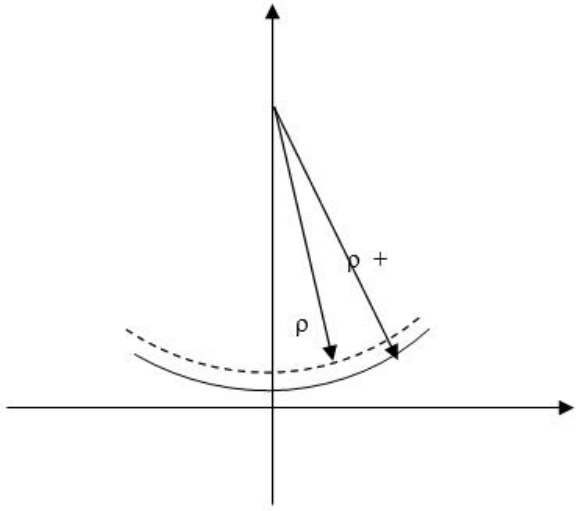}
    \caption{Diameter astrolabe.}
    \label{fig007Diameter}
\end{figure}

However, the centring error can be expressed as an angular error such that 
\begin{equation}~\label{eqq4rate}
dL = \sqrt{p^2d\alpha^2 + dp^2}.
\end{equation}
It is well known that an error in the angle   will cause the arc drawn on the plate to shift, and the resulting error will be non-uniform depending on the radius of the arc and will often be unnoticeable in large arcs. The effect of the error in the radius of the arc will appear large compared to the centring error because it directly affects the determination of the altitude circles, especially in small astrolabes where the spacing is small. If there is an error in    this means that the point falls on another altitude circle (Figure~\ref{fig007Diameter}). For example, in a 15 cm diameter astrolabe, we find that the spacing between the projections of the two arcs $h=0^{\circ}$ and $h=3^{\circ}$ on the midday line is about 2 mm. The radius of the arc $h=0^{\circ}$ is 12.76 cm, meaning that the distance does not exceed 2\% of the radius, which is an amount that could be a drafting error when using traditional methods. This means that if we miscalculate the value   and draw it at a ratio of 2\%, the point under consideration will fall on another arc following the correct arc (Figure~\ref{fig007Diameter}). Consequently, the observation will be on one arc, while the calculation with the astrolabe will be on another adjacent arc, which of course means that the astrolabe is not sufficiently accurate.  
What has been mentioned above applies to any arc of a circle. In the celestial sphere, the error is compounded when determining $p$ due to an error in aligning the east and west lines with the midday line, which causes an error in identifying the center and rotation of the celestial sphere. If an error occurs in the arcs of the hours, it means that there is a compounded error in determining the times of east and west.

\section{Stereographic Design Using Computer Assistance}~\label{sec:Stereographic}

Relying on computer-aided design (CAD) programs for stereographic design has several advantages, as it eliminates some sources of error due to:
\begin{itemize}
    \item Relying on analytical relationships to determine the centers of arcs and the lengths of their diameters,
    \item Drawing based on events and engineering elements for circles and arcs,
    \item Using precise algorithms for dividing arcs and segments into equal parts,
    \item The accuracy of real numbers processed on the computer, which significantly reduces practical error in determining events and measurements.
\end{itemize}

Therefore, using stereographic design programs with computer assistance (CAD System) means eliminating the errors mentioned in points 3 and 4, considering that segmentation operations are carried out using precise algorithms, and the drawing outputs rely on their geometric components. This can be summarized analytically as follows:

The stereographic design relies on two types of engineering procedures:
\begin{itemize}
    \item The main drawing elements (arcs and segments)
    \item Segmentation (division) of engineering components
\end{itemize}

This can be clarified in the following Table~\ref{astdesn1}:
\begin{table}[H]
\centering
\begin{tabular}{|l|l|}
\hline
\multicolumn{2}{|c|}{Astrolabe design}                            \\ \hline
Drawing geometric elements & Segmentation of geometric elements \\ \hline
The three orbits           & Hour curves                        \\ \hline
Arcuate arches             & Arcs of orbits                     \\ \hline
Smoothing brackets         & Zodiac calendar                    \\ \hline
Qibla destination curves   & Sinusoidal descent gradients       \\ \hline
Midday curves              & Square gradients of Shadows        \\ \hline
Recessed brackets          & JAMB gradients   \\ \hline                 
\end{tabular}
\caption{Astrolabe design.}~\label{astdesn1}
\end{table}

\subsection{Drawing the Arcs}

Eliminating the error resulting from determining the center of the arcs is done by calculating its coordinates from the following relations:

\begin{equation}~\label{eqq5rate}
Y_c = \frac{Y_1 + Y_2}{2}
\end{equation}

where

\begin{equation}~\label{eqq6rate}
\left\{
  \begin{array}{l}
    Y_1 = 2R \cot \left( \frac{\varphi + h}{2} \right) \\
    Y_2 = -2R \tan \left( \frac{\varphi - h}{2} \right) \\
  \end{array}
\right.
\end{equation}

\begin{itemize}
  \item $h$: The circle of the required height whose projection is to be drawn.
  \item $\varphi$: The geographical latitude of the country.
  \item $R$: The nominal diameter of the projection (the projection's equator).
\end{itemize}

The coordinates of the assumed point $(X_0, Y_0)$ are set to begin with, and the coordinates of the arc's center are determined by the point $(X_c, Y_c)$.

Eliminating the error resulting from determining the value of the arc's radius can also be done by calculating it from the following relation:

\begin{equation}~\label{eqq7rate}
p=\frac{Y_1-Y_2}{2}.
\end{equation}

\subsection{Drawing the Azimuth Circles}

The first azimuth is drawn, and its center and radius are determined by the following relations:

\begin{equation}~\label{eqq8rate}
Y_c=\frac{Y_1+Y_2}{2},
\end{equation}
where
\begin{equation}~\label{eqq9rate}
\left\{
  \begin{array}{l}
    Y_1=2R \tan \left(45^\circ-\frac{\varphi}{2}\right) \\
    Y_2=-2R \tan \left(45^\circ+\frac{\varphi}{2}\right) \\
  \end{array}
\right.
\end{equation}
and 
\[
p_a=\frac{Y_1+Y_2}{2}.
\]
The remaining smoothing circles have their centers in the same order $Y_c$, and the interval and radius are calculated by the following relations:
\[
\left\{
  \begin{array}{l}
    X_a=p_c\tan A\\
    p=\dfrac{p_c}{\cos A}.
  \end{array}
\right.
\]
Where $A$ is the azimuth value of the circle whose projection is to be calculated from the first azimuth circle.
\subsection{Drawing the curves of the hours}
These curves can be drawn accurately by dividing the arcs of the orbits of Aries, Cancer and Capricorn using division algorithms into twelve equal parts, and by using another algorithm to pass a circular arc from every three points shared in the division order on each of the orbits of Cancer, Equinox and Capricorn. This eliminates the error caused by the points not matching the drawn arc.

\begin{figure}[H]
    \centering
    \includegraphics[width=0.5\linewidth]{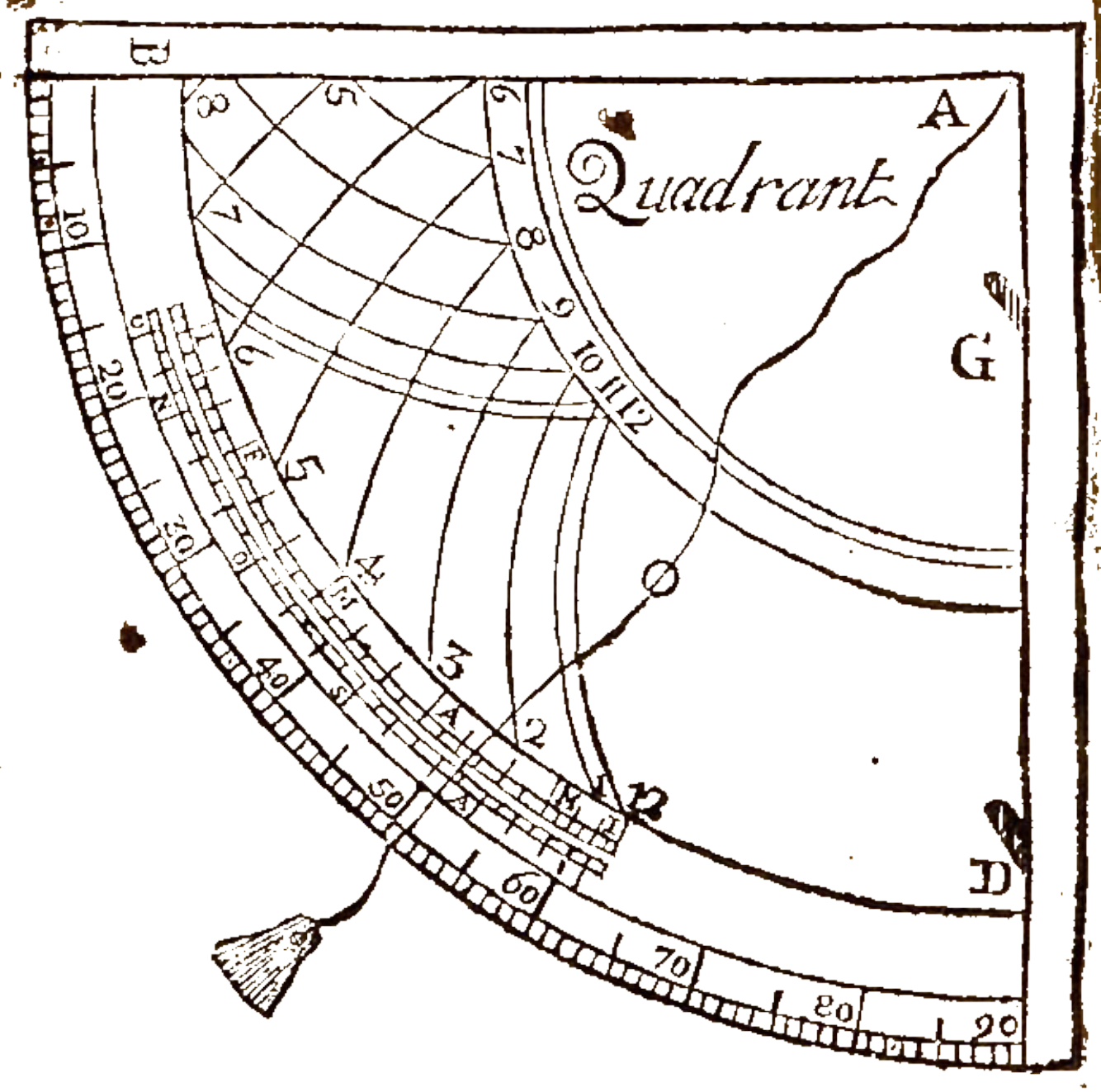}
    \caption{Quadrant.}
    \label{figggquadrantn1}
\end{figure}

\begin{figure}[H]
    \centering
    \includegraphics[width=0.5\linewidth]{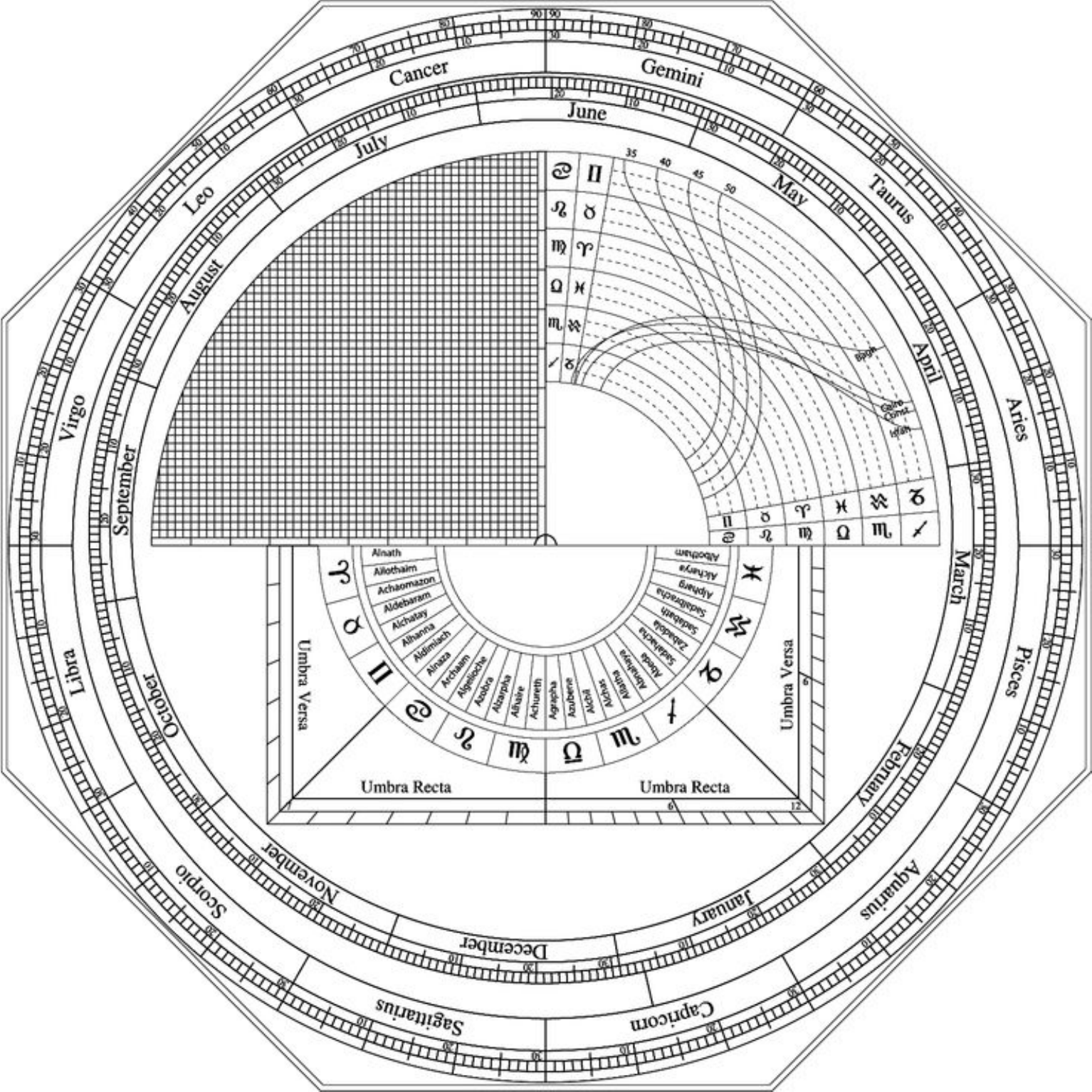}
    \caption{The zodiac calendar.}
    \label{fig00909zodiac}
\end{figure}

\paragraph{Designing the back of the astrolabe:}
 The importance of using computers in the design and manufacture of the back of the astrolabe is evident in the following:
\paragraph{The quadrant:}
The design of the quadrant means dividing its radius into sixty equal parts. In small astrolabes with a diameter of less than 12 cm, the division process becomes difficult using traditional methods because each part is less than a millimetre. Therefore, using a computer and a division algorithm makes the process easy and accurate (Figure~\ref{figggquadrantn1}). 
\paragraph{The zodiac calendar:}
The zodiac calendar is a tool that links the position of the sun in the zodiac circle with time. This means that 365 positions of the sun must be determined according to the days of the year. If we distribute this calendar on the lower half of the astrolabe, this means, according to the previous example, dividing a circle with a length of 18.5 cm into 365 equal parts, which is a very complex operation without a special computer algorithm (Figure~\ref{fig00909zodiac}). In addition, the proguial calendar must be constantly updated because the relationship between the position of the sun and time is not constant due to the precession of the Earth's axis, and this update will only be effective if it is carried out using a computer. 
\subsection{Midday curves and the direction of the Qibla}
\begin{figure}[H]
    \centering
    \includegraphics[width=0.4\linewidth]{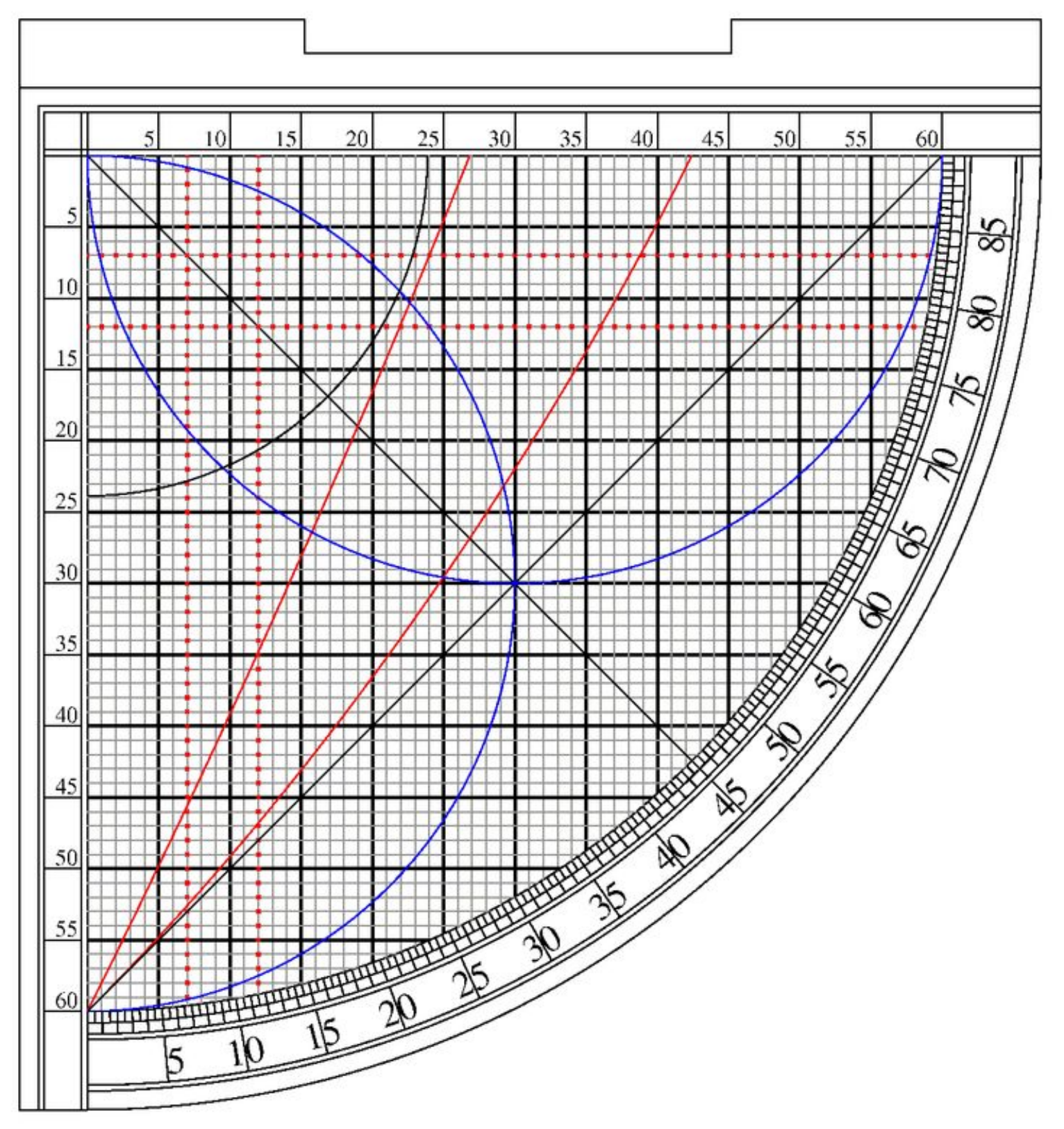}
    \caption{Midday curves}
    \label{figMiddaycurves}
\end{figure}

\paragraph{Midday curves}
The altitude of the sun at midday at latitude $\varphi$ is a function of time, and therefore of the deviation $\delta$, according to the following relationship:
\[
h=90^{\circ}-\varphi+\delta.
\]
The traditional method~\ref{figMiddaycurves} of implementing these curves is to calculate the above equation for three values of deviation $\delta$, and then search for the best arc centre passing through the three points. This process will of course be easier, faster and more accurate using one of the $ACAD$ software algorithms.

\subsection{Experiments and the Validity of the Model}

These experiments aim to verify the relationship between the angle of inclination of the tube and the contact angle. The following mathematical relationships are used:

\begin{equation}~\label{eq1Validity}
\tan \alpha = \frac{\cos \varphi_M \sin (\lambda_M - \lambda)}{\cos \varphi_M \sin \varphi - \sin \varphi_M \cos \varphi \cos (\lambda_M - \lambda)} 
\end{equation}

\begin{equation}~\label{eq2Validity}
\sin \delta = \sin \varphi \sin h + \cos \varphi \cosh \cos \alpha \tag{14}
\end{equation}

\noindent where:
\begin{itemize}
    \item $\varphi, \lambda$: coordinated with direction of sun.
    \item $\varphi_M, \lambda_M$: coordinates of the observation point by considering Mecca.
\end{itemize}

The experiment was conducted in two stages. First, the capillary tube was kept vertical to measure the initial contact angle for the liquid (water or mercury). Then, using \texttt{MATLAB}, the following were done:

\begin{itemize}
    \item Data from the measurements were entered to create a program that calculates the contact angle using the above equations.
    \item The tube was tilted at specific angles, and the change in the contact angle was observed.
    \item The results showed that the contact angle \emph{increases} slightly with increasing tilt angle, but this increase is small and can be \emph{neglected} in most calculations.
\end{itemize}

\section{Conclusion}
Through this paper, the astrolabe is considered one of the most important astronomical instruments that contributed to the development of astronomy, as it is regarded as the foundation upon which astronomy was built and served as a kind of computer that made life easier at that time, which was the golden age of the astrolabe.  We presented an analytical study of the method of making the astrolabe by engineering and software methods due to the historical importance through which the astrolabe contributed to the development of astronomy over several time periods.

\end{document}